\documentclass[12pt]{amsart}
\newcommand{\abs}[1]{\ensuremath{\left|{#1}\right|}}

\newcommand{\mg}{\infty}

\newcommand{\pii}{\pi i}

\newcommand{\piz}{\pi z}

\newcommand{\qte}[1]{\q\te{#1}}

\renewcommand{\b}{\mathbf{b}}

\newcommand{\dsum}{\di\sum}

\newcommand{\hgt}[1]{\rule[0pt]{0pt}{#1}}

\newcommand{\dep}[1]{{\vrule width 0pt height 0pt depth #1}}

\DeclareFontFamily{U}{mathx}{\hyphenchar\font45}
\DeclareFontShape{U}{mathx}{m}{n}{
      <5> <6> <7> <8> <9> <10>
      <10.95> <12> <14.4> <17.28> <20.74> <24.88>
      mathx10
      }{}
\DeclareSymbolFont{mathx}{U}{mathx}{m}{n}
\DeclareFontSubstitution{U}{mathx}{m}{n}
\DeclareMathAccent{\widecheck}{0}{mathx}{"71}

\newcommand{\inrr}{\ensuremath{\in\rr}}

\renewcommand{\kill}[1]{}
\newcommand{\dummy}[1]{\mbox{}}

\makeatletter
\newcommand{\xequal}[2][]{\ext@arrow 0055{\equalfill@}{#1}{#2}}
\def\equalfill@{\arrowfill@\Relbar\Relbar\Relbar}
\makeatother

\newcommand{\xeq}{\xequal}

\newcommand{\mto}{\mapsto}

\newcommand{\Set}[2]{\ensuremath{\left\{{#1}\,\middle|\,{#2}\right\}}}

\renewcommand{\k}{\ensuremath{\ol{\mathrm{P}}}}

\newcommand{\h}{\hline}

\renewcommand{\k}[1]{\ensuremath{\left({#1}\right)}}

\newcommand{\re}{\te{Re}\,}

\newcommand{\bca}{\begin{cases}}
\newcommand{\eca}{\end{cases}}

\newcommand{\mug}{\ensuremath{\infty}}

\newcommand{\ff}[2]{\ensuremath{\di\fr{#1}{#2}}}

\newcommand{\bpic}{\begin{picture}}\newcommand{\epic}{\end{picture}}

\newcommand{\beda}{\begin{edaenumerate}}
\newcommand{\eeda}{\end{edaenumerate}}

\newcommand{\g}{\ensuremath{\mathbf{g}}}

\newcommand{\cd}{\cdots}

\newcommand{\q}{\quad}

\newcommand{\bq}{\begin{quote}}\newcommand{\eq}{\end{quote}}

\newcommand{\rt}{\sqrt}
\newcommand{\be}{\begin{enumerate}}\newcommand{\ee}{\end{enumerate}}
\newcommand{\bce}{\begin{center}}\newcommand{\ece}{\end{center}}
\newcommand{\bde}{\begin{description}}\newcommand{\ede}{\end{description}}
\newcommand{\bri}{\begin{flushright}}\newcommand{\eri}{\end{flushright}}
\newcommand{\bb}{\begin{block}}\newcommand{\eb}{\end{block}}
\newcommand{\bt}{\begin{thm}}\newcommand{\et}{\end{thm}}
\newcommand{\bpf}{\begin{proof}}\newcommand{\epf}{\end{proof}}
\newcommand{\bex}{\begin{ex}}\newcommand{\eex}{\end{ex}}
\newcommand{\bexr}{\begin{exr}}\newcommand{\eexr}{\end{exr}}
\newcommand{\bft}{\begin{fact}}\newcommand{\eft}{\end{fact}}
\newcommand{\brk}{\begin{rmk}}\newcommand{\erk}{\end{rmk}}
\newcommand{\ba}{\begin{align*}}\newcommand{\ea}{\end{align*}}
\newcommand{\bexe}{\begin{exe}}\newcommand{\eexe}{\end{exe}}

\newcommand{\bit}{\begin{itemize}}\newcommand{\eit}{\end{itemize}}

\newcommand{\bcm}{}

\newcommand{\hf}{\hfill}
\newcommand{\ci}{\CIRCLE}
\newcommand{\fr}{\frac}
\newcommand{\cc}{\ensuremath{\mathbf{C}}}

\newcommand{\rr}{\ensuremath{\mathbf{R}}}

\newcommand{\bd}{\begin{defn}}\newcommand{\ed}{\end{defn}}
\newcommand{\bp}{\begin{prop}}\newcommand{\ep}{\end{prop}}
\newcommand{\p}{\ensuremath{\pi}}
\newcommand{\eh}{\emph}\newcommand{\al}{\alpha}

\newcommand{\te}{\text}\newcommand{\ph}{\phantom}

\newcommand{\di}{\displaystyle}\renewcommand{\a}{\ensuremath{\bm{a}}}
\renewcommand{\b}{\ensuremath{\bm{b}}}

\newcommand{\f}{\frac}

\newcommand{\z}{\ensuremath{\bm{z}}}
\newcommand{\np}{\newpage}

\renewcommand{\b}{\beta}

\renewcommand{\a}{\alpha}

\renewcommand{\P}{\mathcal{P}}

\usepackage[dvipdfmx]{graphicx}
\usepackage[dvipsnames]{xcolor}
\usepackage{float,afterpage}
\usepackage{asymptote,layout}
\usepackage{wrapfig,epic}
\usepackage{amsmath,amssymb,cases,pict2e}

\usepackage{multicol}
\usepackage{tikz}

\usepackage{geometry}
\usepackage{exscale,latexsym,bm}
\usepackage{amssymb,enumerate,amsmath,amsthm,amsfonts}
\usepackage{verbatim,fancybox,wasysym,fancyhdr,type1cm}
\usepackage[frame,all,poly,curve,knot,arrow]{xy}
\usepackage{colortbl}
\usepackage{boxedminipage}
\usepackage{multirow}

\graphicspath{{./figs/}}
\theoremstyle{definition}
\newtheorem{thm}{Theorem}[section]

\newtheorem{lem}[thm]{Lemma}
\newtheorem{prop}[thm]{Proposition}

\newtheorem{exr}[thm]{Exercise}

\newtheorem{ex}[thm]{Example}

\newtheorem{defn}[thm]{Definition}\newtheorem{rmk}[thm]{Remark}
\newtheorem{fact}[thm]{Fact}
\newtheorem{block}[thm]{}
\newtheorem*{exe}{Exercise}

\renewcommand{\a}{\alpha}

\renewcommand{\h}{\hline}
\renewcommand{\arraystretch}{1.5}

\begin{document}

\title{Values of zeta-one functions at positive even integers}
\author{Masato Kobayashi}
\author{Shunji Sasaki}
\thanks{corresponding author:Masato Kobayashi, masato210@gmail.com}
\date{\today}                                       
\address{Masato Kobayashi\\
Department of Engineering\\
Kanagawa University, 3-27-1 Rokkaku-bashi, Yokohama 221-8686, Japan.}
\email{masato210@gmail.com}
\address{Shunji Sasaki\\
Kawaguchi public Kamiaoki junior high school\\
3-9-1 Kamiaoki-Nishi, Kawaguchi 333-0845, Japan.
}

\email{schnittkejp@me.com}

\maketitle
\begin{abstract}
Motivated by Euler-Goldbach and Shallit-Zikan theorems, we introduce zeta-one functions with infinite sums of  $n^{s}\pm1$ as an analogy of the Riemann zeta  function. Then we compute values of these functions at positive even integers by the residue theorem.
\end{abstract}
\tableofcontents
\np

\renewcommand{\g}{\gamma}
\newcommand{\Om}{\Omega}
\renewcommand{\P}{\mathbf{P}}
\renewcommand{\z}{\zeta}

\newcommand{\vp}{\zeta_{1}}
\renewcommand{\re}{\text{Re}}
\newcommand{\tf}{\tfrac}
\renewcommand{\g}{\Gamma}
\newcommand{\res}{\text{Res}}

\section{Introduction: Euler-Goldbach Theorem}
Let us start with the celebrated \eh{Euler-Goldbach Theorem}. Say that a natural number $p$ is a \eh{perfect power} if $p=n^{m}$ for some natural numbers $m, n\ge2$.

\begin{thm}[Euler-Goldbach]
\[
\dsum_{p:\text{perfect power}}\ff{1}{p-1}=1.
\]
\end{thm}
See Bibiloni-Paradis-Viader \cite{bibiloni} for history of  this theorem.
Recently, Shallit-Zikan \cite{shallit} (1983) reinterpreted it in terms of \eh{Riemann's zeta function}:
the infinite series 
\begin{align*}
	\zeta(s)&=\dsum_{n=1}^{\mug}\ff{1}{n^{s}}=
	\ff{1}{1^{s}}+\ff{1}{2^{s}}+\ff{1}{3^{s}}+\ff{1}{4^{s}}+\ff{1}{5^{s}}+\cd
\end{align*}
is convergent for all complex numbers $s$ such that $\re{(s)}>1$. 
Indeed, Euler proved that 
\[
\z(2k)=-\ff{1}{2}\ff{(2\pii)^{2k}}{(2k)!}B_{2k}
\q\te{(Table \ref{zeven})}
\]
where $\{B_{n}\}$ are signed Bernoulli numbers as in 
Table \ref{bno}; refer to Ayoub \cite{ayoub} for history of this function.
\newcommand{\tff}{\tfrac}

Since $\zeta(s)=1+\tff{1}{2^{s}}+\cd>1$ and 
\[
2>
\ff{\p^{2}}{6}=\z(2)>\z(3)>\z(4)>\z(5)>\cd,
\]
we have $1<\z(s)<2$ for all $s\ge 2$.
That is, $\z(s)-1$ is the \eh{fractional part} of $\z(s)$. For example, 
\begin{align*}
	\zeta (2)-1&=0.6449\dots,
	\\\zeta (3)-1&=0.2020\dots,
	\\\zeta (4)-1&=0.0823\dots,
	\\\zeta (5)-1&=0.0369\dots.
\end{align*}
\begin{thm}[Shallit-Zikan \cite{shallit}]
\begin{align*}
	\dsum_{k=2}^{\mug}(\z(k)-1)&=1.
\end{align*}
\end{thm}

Let us now see a similar result.


\begin{thm}[{See J.M. Borwein-Bradley-Crandall \cite[p. 262]{bbc}}]
\[
\dsum_{k=1}^{\mug}(\z(2k)-1)=\ff{3}{4}.
\]
\end{thm}
Here we give a proof since it suggests some ideas for our main results.
\begin{proof}
Consider the double sequence 
$a_{nk}=(\tff{1}{n^{2k}})_{n\ge2, k\ge1}$
and positive series $\sum_{n\ge2, k\ge1}a_{nk}$.
We find that 
\begin{align*}
	\dsum_{n=2}^{\mug} \dsum_{k=1}^{\mug}\ff{1}{n^{2k}}
&=
	\dsum_{n=2}^{\mug} \ff{1}{n^{2}}
	\dsum_{k=0}^{\mug}\k{\ff{1}{n^{2}}}^{k}
	\\&=\dsum_{n=2}^{\mug} \ff{1}{n^{2}}
	\ff{1}{1-\f{1}{n^{2}}}
	= \dsum_{n=2}^{\mug} \ff{1}{n^{2}-1}
	\\&=
	\ff{1}{2}
	\dsum_{n=2}^{\mug} \k{\ff{1}{n-1}-\ff{1}{n+1}}=
	\ff{3}{4}
\end{align*}
 so that we can freely switch order of this series.
As a consequence, 
\begin{align*}
	\dsum_{k=1}^{\mug}(\z(2k)-1)&=
	\dsum_{k=1}^{\mug} \dsum_{n=2}^{\mug} \ff{1}{n^{2k}}=
	\dsum_{n=2}^{\mug} \dsum_{k=1}^{\mug}\ff{1}{n^{2k}}=\ff{3}{4}.
	\end{align*}
\end{proof}
In this proof, the infinite series 
$\sum_{n=2}^{\mug} \tf{1}{n^{2}-1}$ appeared.
It is now natural to think of analogous sums 
$\sum\tf{1}{n^{s}\pm1}$.
With this simple idea, this article introduces \eh{zeta-one functions} $\z_{+1}(s)$, $\z_{-1}(s)$ and we compute values of $\z_{+1}(2m)$ and $\z_{-1}(2m)$ as main Theorems \ref{th1} and \ref{th2}.

%



\section{Zeta-one functions}
\subsection{Definition}

Throughout $k, m, n, N$ and $s$ each denote a nonnegative integer unless otherwise specified. Further, we assume that $s\ge2$.
\begin{defn}
Define the \eh{zeta-one functions} by 
\[
\z_{+1}(s)=\dsum_{n=1}^{\mug}\ff{1}{n^{s}+1} 
\qte{and}\q
\z_{-1}(s)= \dsum_{n=2}^{\mug} \ff{1}{n^{s}-1}
\]
(For $s\ge2$, these sums are indeed convergent as mentioned below). Call each \eh{zeta-plus-one} and \eh{zeta-minus-one} function, respectively. 
\end{defn}

\subsection{Example}

\begin{ex}
As seen above, $\z_{-1}(2)=\tff{3}{4}$. Moreover, 
since 
\[
\coth(\piz)=\ff{1}{\piz}+\ff{1}{\p}
\dsum_{n=1}^{\mug}\ff{2z}{z^{2}+n^{2}},
\q z\in\cc,
\]
the substitution $z=1$ implies that 
\[
\z_{+1}(2)=
\dsum_{n=1}^{\mug}\ff{1}{n^{2}+1}=
-\ff{1}{2}+\ff{\p}{2}\coth(\p).
\]
As a consequence, 
\[
\z_{+1}(s)\le \z_{+1}(2)<\mug, \q 
\z_{-1}(s)\le \z_{-1}(2)<\mug
\]
for all $s\ge2$.
\end{ex}


{\renewcommand{\arraystretch}{2.15}
\begin{table}
\caption{zeta even values}
\index{Z2n@$\zeta(2s)$}
\label{zeven}
\begin{center}
\begin{tabular}{c|ccccccccccccccccccc}
$2k$&2&4&6&8&10&12&
$\cd$\\\hline
$\z(2k)$&
$\ff{\p^{2}}{6}$&
$\ff{\p^{4}}{90}$&
$\ff{\p^{6}}{945}$&
$\ff{\p^{8}}{9450}$&
$\ff{\p^{10}}{93555}$&
$\ff{691\p^{12}}{638512875}$&
$\cd$\\
\end{tabular}
\end{center}
\end{table}%
}
{\renewcommand{\arraystretch}{2}
\begin{table}
\caption{signed Bernoulli numbers 
$(\te{with } B_{3}=B_{5}=\cd =0$)}
\label{bno}
\begin{center}
\begin{tabular}{c|ccccccccccccccccccccccccc}
$n$&0&1&2&4&6&8\\\h
$B_{n}$&\dep{17pt}1&$- \f12$&$\f{1}{6}$\hgt{20pt}&$-\f{1}{30}$
&$\f{1}{42}$&$-\f{1}{30}$\\\hline
$n$&10&12&14&16&18&20\\\hline
$B_{n}$&\dep{17pt}$\f{5}{66}$&$-\f{691}{2730}$&$\f{7}{6}$&$-\f{3617}{510}$&
$\f{43867}{798}$
&$-\f{174611}{330}$\\\h
\end{tabular}
\end{center}
\end{table}%
}

\section{Main theorem 1}

Toward the proof of Theorem \ref{th1} on $\z_{+1}(2m)$, we need lemmas.

\subsection{Lemmas}

For $m\ge1$, set
\[
f(z)=\frac{\cot(\pi z)}{z^{2m}+1}
\qte{and}\q \al=\exp\k{\ff{\pi i}{2m}}.
\]
Recall that 
\[
 \cot(\pi z)=\frac{1}{\pi z}+ \frac{1}{\pi}\sum_{n=1}^\infty \frac{2z}{z^2-n^2}, 
 \q z\in \cc.
 \]
Thus, $f$ is a meromorphic function with simple poles $z=0,\pm1,\pm2,\pm3,\dotsb$ and 
\begin{align*}
\a,\a^3,\a^5, \dotsb, \a^{2m-1},-\a,-\a^3,-\a^5, \dotsb,  -\a^{2m-1}
 \end{align*}
as all the roots of $z^{2m}+1=0$. Let us compute the residue of $f$ at each pole.

\begin{lem}\label{l2}
For $n=0, \pm1, \pm2, \pm3, \dotsb$, 
we have 
%
\[
\text{Res}(f, n)=\ff{1}{\p(n^{2m}+1)},
\]
and for $1\le k\le m$,
\[
\text{Res}(f,\pm\al^{2k-1})=
-
\ff{\a^{2k-1}{\cot(\pi \a^{2k-1})}}{2m}.
\]
\end{lem}
\begin{proof}
First, we have 
\begin{align*}
 \text{Res}(f,n)=&\lim_{z\to n}(z-n)f(z){}=
 \lim_{z\to n}\ff{z-n}{\sin\pi (z-n)}\cdot\ff{\cos(\piz)}{z^{2m}+1}{}
 =\ff{1}{\pi(n^{2m}+1)}.
\end{align*}
%
Second, for $1\le k\le m$,
\begin{align*}
	 \text{Res}(f,\al^{2k-1})&=
	 \lim_{z\to \a^{2k-1}}(z-\a^{2k-1})f(z)
	\\&=
	 \lim_{z\to \a^{2k-1}}(z-\a^{2k-1})\ff{\cot(\pi z)}{
	 z^{2m}+1}	
	\\&=
	\lim_{z\to \a^{2k-1}}\cot(\pi z)
	\lim_{z\to \a^{2k-1}}\ff{z-\a^{2k-1}}{
	 z^{2m}+1}
	\\&=
	{\cot(\pi \a^{2k-1})}
	\lim_{z\to \a^{2k-1}}\ff{1}{
	 2mz^{2m-1}} \q \te{(L'H\^{o}pital's rule)}
	\\&=
		{\cot(\pi \a^{2k-1})}
	\ff{1}{
	 2m(\al^{2k-1})^{2m-1}}
	\\&=-
\ff{\a^{2k-1}{\cot(\pi \a^{2k-1})}}{2m} \q{\k{(\a^{2k-1})^{2m}=-1}
}.
\end{align*}

%
It is quite similar to show that 
\[
\te{Res}(f, -\a^{2k-1})=-\ff{\a^{2k-1}\cot(\p\a^{2k-1})}{2m}.
\]

\end{proof}

\begin{lem}\label{l1}
For a positive integer $N$, consider line segments on the complex plane 
\begin{align*}
	C_{1}(N)&=
	\Set{\k{N+\f{1}{2}}+yi
}{-\k{N+\f{1}{2}}\le y\le {N+\f{1}{2}}},
	\\C_{2}(N)&=
	\Set{x+\k{N+\f{1}{2}}i
}{-\k{N+\f{1}{2}}\le x\le {N+\f{1}{2}}},
	\\C_{3}(N)&=
	\Set{-\k{N+\f{1}{2}}+yi
}{-\k{N+\f{1}{2}}\le y\le {N+\f{1}{2}}},
	\\C_{4}(N)&=
	\Set{x-\k{N+\f{1}{2}}i
}{-\k{N+\f{1}{2}}\le x\le {N+\f{1}{2}}}
\end{align*}
and set $C(N)=C_{1}(N)\cup C_{2}(N)\cup C_{3}(N)\cup 
C_{4}(N).$
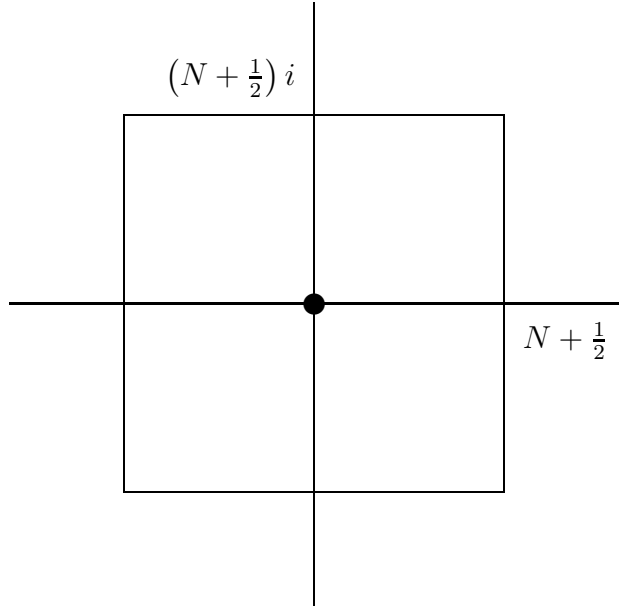
\[
\begin{xy}
0;<5mm,0mm>:
(5,0)*{}="al1"*++!UL{N+\f{1}{2}}
,(0,5)*{}="al2"*++!DR{\k{N+\f{1}{2}}i
}
,(0,0)*{\ci}="t1"*+++!D{}
,\ar@{-}(-8,0);(8,0)
,\ar@{-}(0,-8);(0,8)
,\ar@{-}(-5,5);(5,5)
,\ar@{-}(5,5);(5,-5)
,\ar@{-}(5,-5);(-5,-5)
,\ar@{-}(-5,-5);(-5,5)
\end{xy}\]
\begin{enumerate}
	\item If $z\in C(N)$, then 
	$|\cot(\piz)|\le \coth{\tff{3}{2}\p}$.
	\item If $z\in C(N)$, then 
	\[
\ff{1}{\abs{z^{2m}+1}}\le \ff{1}{
	\k{N+\f{1}{2}}^{2m}-1}.
	\]
\end{enumerate}
\end{lem}

\begin{proof}\hf
\begin{enumerate}
	\item Suppose $z\in C(N)$.
	If $z\in C_{1}(N)$, then write 
	\[
z=\k{N+\f{1}{2}}+yi, \q -\k{N+\f{1}{2}}\le y\le {N+\f{1}{2}}.\]
\begin{align*}
	\abs{\cot\piz}&=
	\abs{
	\ff{e^{\pii z}+e^{-\pii z}
	}{
	e^{\pii z}-e^{-\pii z}}
	}
	=
	\abs{
	\ff{
	e^{-\pi y}e^{\k{N+\f{1}{2}}\pii}
	+e^{\pi y}e^{-\k{N+\f{1}{2}}\pii
	}}{
	e^{-\pi y}e^{\k{N+\f{1}{2}}\pii}
	-e^{\pi y}e^{-\k{N+\f{1}{2}}\pii
	}
	}
	}
	\\&=
	\abs{
	\ff{
	e^{-\pi y}(-1)^{N}i
	+e^{\pi y}(-1)^{N}(-i)
	}
	{
	e^{-\pi y}(-1)^{N}i
	-e^{\pi y}(-1)^{N}(-i)
	}
	}
	\\&=
	\abs{
	\ff{e^{-\pi y}-e^{\pi y}
	}
	{
	e^{-\pi y}+e^{\pi y}}
	}
	=
	\abs{
	\ff{e^{\pi y}-e^{-\pi y}
	}
	{
	e^{\pi y}+e^{-\pi y}
	}
	}
	\\&=|\tanh(y)|\le 1<
	\coth\ff{3}{2}\p\,\, (=1.00016\cd).
\end{align*}
If $z\in C_{2}(N)$, then 
\[
z=x+
\k{N+\f{1}{2}}i, \q -\k{N+\f{1}{2}}\le x\le {N+\f{1}{2}}
\]
and \begin{align*}
	|\cot(\piz)|&=
	\abs{
	\ff{e^{\pii z}+e^{-\pii z}}
	{e^{\pii z}-e^{-\pii z}}
	}
	\\&\le 
	\ff{\abs{e^{\pii z}}+\abs{e^{-\pii z}}}
	{
\abs{
\abs{e^{\pii z}}
-
\abs{e^{-\pii z}}
}
}
	\\&=
	\ff{
	\abs{
	e^{-\pi \k{N+\f{1}{2}}}}
	+\abs{
	e^{\pi \k{N+\f{1}{2}}}
	}
	}
	{
	\abs
	{\abs{
	e^{-\pi \k{N+\f{1}{2}}}
	}
	-
	\abs
	{e^{\pi \k{N+\f{1}{2}}}
	}}
	}
	\\&=\ff{e^{\pi \k{N+\f{1}{2}}}+e^{-\pi \k{N+\f{1}{2}}}
	}{e^{\pi \k{N+\f{1}{2}}}-e^{-\pi \k{N+\f{1}{2}}}}
	\\&=\coth\k{N+\f{1}{2}}\p\le \coth\ff{3}{2}\p
\end{align*}
since $t\mto \coth(t)$ is decreasing for $t>0$.
	For $z\in C_{3}(N)\cup C_{4}(N)$, 
	we have  $-z\in C_{1}(N)\cup C_{2}(N)$ so that 
\[
	|\cot(\piz)|=
	|-\cot(\pi (-z))|=
	|\cot(\pi (-z))|\le \coth\ff{3}{2}\p.
\]
	\item If $z\in C(N)$, then 
	$|z|\ge N+\f{1}{2}$.
	Consequently, 
$|z|^{2m}-1\ge \k{N+\f{1}{2}}^{2m}-1
$,
\[
\ff{1}{|z|^{2m}-1}\le \ff{1}{\k{N+\f{1}{2}}^{2m}-1},
\]
\[
\ff{1}{|z^{2m}+1|}\le 
\ff{1}{|z|^{2m}-1}\le \ff{1}{\k{N+\f{1}{2}}^{2m}-1}.
\]
\end{enumerate}

\end{proof}

\subsection{Proof of main theorem 1}

Let $s$ be a positive even integer, say $s=2m$, $m\ge1$.
Further, let $\al=\al_{2m}=\exp\k{\ff{\pii}{2m}}$ 
for convenience.
\begin{thm}
\label{th1}
\begin{align*}
\z_{+1}(2m)&=
-\frac{1}{2}+\frac{1}{2m}\sum_{k=1}^{m}
\pi\al^{2k-1}\cot(\pi \a^{2k-1}).
 \end{align*}
\end{thm}
\begin{proof}
View $C(N)=C_{1}(N)+C_{2}(N)+C_{3}(N)+C_{4}(N)$ above as the sum of four paths with counterclockwise orientation. 
Notice that any pole of $f$ does not lie on $C(N)$. 
We are going to compute the integral 
\[
I_{N}=\di\int_{C(N)}f(z)dz.
\]
Let $D(N)$ be the domain enclosed by $C(N)$. Then, 
the residue theorem with Lemma \ref{l2} claims that 
\begin{align*}
	I_{N}&=
	\di\int_{C{(N)}}f(z)dz
	\\&=
	2\pii
	\sum_{\substack{a: \text{pole of\ph{a}} f(z)\\
	a\in D(N)}
	}\text{Res}(f, a)
	\\&=
	2\pii
	\k{\text{Res}(f, 0)+
	\dsum_{n=1}^{N}
	\k{\text{Res}(f, n)+\text{Res}(f, -n)}
	+
	\dsum_{k=1}^{m}
	\k{\text{Res}(f, \a^{2k-1})+\text{Res}(f, -\al^{2k-1})}}
	\\&=
	2\pii
\k{\ff{1}{\p}+2
\dsum_{n=1}^{N}\ff{1}{\p(n^{2m}+1)}
+2\k{
-
\dsum_{k=1}^{m}\ff{\a^{2k-1}\cot(\pi\a^{2k-1})}{2m}
}}
\end{align*}
while Lemma \ref{l1} implies 
\begin{align*}
	|I_{N}|&=\abs{\di\int_{C{(N)}}f(z)dz}
	\\&\le \di\int_{C{(N)}}\abs{f(z)}dz
	\\&\le 
	\ff{\coth\f{3}{2}\p}{
\k{N+\f{1}{2}}^{2m}-1
}
\di\int_{C{(N)}}dz
	\\&=\ff{\coth\f{3}{2}\p}{
\k{N+\f{1}{2}}^{2m}-1
}\cdot 8\k{N+\f{1}{2}}\to 0 \q (N\to \mg).
\end{align*}
%
Therefore, taking the limit $N\to\mg$ for $I_{N}$ yields 
\[
0=
2\pii
\k{\ff{1}{\p}+\ff{2}{\p}\z_{+1}(2m)-2
\dsum_{k=1}^{m}\ff{\a^{2k-1}\cot(\pi\a^{2k-1})}{2m}}
.
\]
Conclude that 
\[
\z_{+1}(2m)=
-\frac{1}{2}+\frac{1}{2m}\sum_{k=1}^{m}
\pi\al^{2k-1}\cot(\pi \a^{2k-1}).
\]

\end{proof}

\subsection{Example}

Of course, $\z_{+1}(2m)$ is a real number so that there should be some expression of $\z_{+1}(2m)$ in terms of only real trigonometric functions.
\begin{ex}
Let $s=4, m=2$ and $\a=\a_{4}=\exp\k{\tff{\pi i}{4}}$. Then 
   \begin{align*}
  \z_{+1}(4)=\sum_{n=1}^\infty \frac{1}{n^{4}+1}&=-\frac{1}{2}+\frac{\pi}{4}\sum_{k=1}^2 
  \a^{}\cot(\pi \a^{}){}\\
  &=-\frac{1}{2}+\frac{\pi}{4}
   \k{\a \cot(\pi \a^{})+\a^{3} \cot(\pi \a^{3})}
   {}\\
  &=-\frac{1}{2}+\frac{\pi}{4}
 \k{\a \cot(\pi \a)+\a^{-1} \cot(\pi \a^{-1})}
.
   \end{align*}
Now, it follows from the facts 
\[\pi \exp\k{\pm\frac{\pi i}{4}}  =\frac{1\pm i}{\sqrt2} \pi
\q \te{and} \q \cot(x+yi)=\ff{\sin2x-i\sinh2y}{\cosh2y-\cos2x} 
\q x, y\inrr
\]
that
 \begin{align*}
 \z_{+1}(4)&=
 	-\frac{1}{2}+\frac{\pi}{4}
 \k{\a \cot(\pi \a)+\a^{-1} \cot(\pi \a^{-1})}
  \\&=
 -\frac{1}{2}+
 \frac{\pi}{4}
\k{
\ff{1+i}{\rt{2}}\cdot\ff{\sin\rt{2}\p-i\sinh\rt{2}\p}{\cosh\rt{2}\p-\cos\rt{2}\p}
+
\ff{1-i}{\rt{2}}\cdot\ff{\sin\rt{2}\p+i\sinh\rt{2}\p}{\cosh\rt{2}\p-\cos\rt{2}\p}
}
 	\\&=
	-\frac{1}{2}+\frac{\sqrt2 \pi}{4}
  \k{\frac{ \sin\sqrt2 \pi+ \sinh\sqrt2 \pi}{\cosh\sqrt2 \pi-\cos\sqrt2 \pi}}.
 \end{align*}

Let $s=6$, $m=3$ and 
$\a=\a_{6}=\exp\k{\tff{\pii}{6}}=\tff{\rt{3}+i}{2}$. 
With $\a^{3}=i$, $\a^{5}=\a^{-1}=\tff{\rt{3}-i}{2}$ and 
$i\cot(\pi i)=\coth(\p) $, we observe that 
\begin{align*}
  \z_{+1}(6)&=
  -\frac{1}{2}+\frac{\pi}{6}
  \k{\al\cot(\pi \a)+
  {\al^{3}\cot(\pi \a^{3})}
+\al^{5}\cot(\pi \a^{5})}
\\
  &=
  -\frac{1}{2}+\frac{\pi}{6}
  \k{\al\cot(\pi \a)+
  \al^{-1}\cot(\pi \a^{-1})
+ i\cot(\pii)
}
\\
&= -\frac{1}{2}+\frac{\pi}{6}
\k{
\ff{\rt{3}+i}{2}\cdot
\ff{\sin\rt{3}\p-i\sinh\p}{\cosh\p-\cos\rt{3}\p}
+
\ff{\rt{3}-i}{2}\cdot
\ff{\sin\rt{3}\p+i\sinh\p}{\cosh\p-\cos\rt{3}\p}
+\coth(\p)
}
\\&=
-\frac{1}{2}+\frac{\pi}{6}\k{
\frac{\sqrt3 \sin\sqrt3 \pi+\sinh\pi}{\cosh\pi-\cos\sqrt3 \pi }+\coth(\pi)}.
 \end{align*}
%


   \end{ex}

\section{Main theorem 2}

Next, we prove Theorem \ref{th2} on $\z_{-1}(2m)$. 
Ideas are quite same.

\subsection{Lemma}

For $m\ge1$, let 
\[
g(z)=\ff{\cot(\piz)}{z^{2m}-1} \qte{and}\q 
\b=\exp\k{\ff{\pii}{m}}.
\]

It has poles at 
$z=0, \pm1, \pm 2, \pm3, \dots,$ 
and $z=\b^{k}$, $1\le k\le 2m-1, k\ne m$. 
The order of the poles $z=\pm1$ is \eh{2} and the all others are simple.

%

\begin{lem}\label{l3}
For $n=0, \pm2, \pm3, \dots$, we have 
\[
\res(g, n)=\ff{1}{\p(n^{2m}-1)}, 
\]
for $1\le k\le 2m-1, k\ne m$, 
%
\[\res(\p, \b^{k})=\ff{\b^{k}}{2m}\cot(\pi \b^{k})
\]
and moreover 
\[
\res(g, \pm1)=-\ff{2m-1}{4m\p}.
\]
\end{lem}

\begin{proof}
The proofs of the first two equalities 
are almost similar to ones for Lemma \ref{l2}.
Thus we only need to show 
$\res(g, \pm1)=-\tff{2m-1}{4m\p}$.

Let $\phi(z)= \sum_{k=0}^{2m-1}z^{k}$. 
Notice that $\phi(z)=(z^{2m+1}-1)/(z-1)$.

Then
\begin{align*}
	\res(g, 1)&=\di\lim_{z\to 1}{ \ff{d}{dz}(z-1)^{2}g(z)}
	\\&=\di\lim_{z\to 1}{ \ff{d}{dz}(z-1)\cot(\p z)\cdot\ff{1}{\phi(z)}}
	\\&=
	\di\lim_{z\to 1}\k{ 
	\ff{\cot(\pi z)-\p(z-1)(\cot^{2}\pi z+1)}{\phi(z)}
	-
	(z-1)\cot(\pi z)\ff{\phi'(z)}{\phi(z)^{2}}
	}.
\end{align*}

Let us see the first term. 
Immediately, $\lim_{z\to 1}\phi(1)=2m$ and 
\begin{align*}
	\di\lim_{z\to 1}
	\k{
\cot(\pi z)-\p(z-1)(\cot^{2}\pi z+1)}
&
\xeq{w=z-1}
\di\lim_{w\to 0}\k{
\cot(\pi w)-\pi w(\cot^{2}\pi w+1)}
	\\&=
	\di\lim_{w\to 0}
	\k{
\cot(\pi w)(1-\pi w\cot(\pi w))}
-\di\lim_{w\to 0}{\pi w}
	\\&=
	\di\lim_{w\to 0}
	{
\ff{\tan(\pi w)-\pi w}{\tan^{2}\pi w}}
-0
	\\&=
	\di\lim_{w\to 0}{
\ff{\pi(1+\tan^{2}\pi w)-\pi}
{2\pi \tan \pi w(1+\tan^{2}\pi w)}
}
=0. \q \te{(L'H\^{o}pital's rule)}
\end{align*}
In addition, since 
\[
\phi'(1)= \dsum_{k=0}^{2m-1}k
=m(2m-1),
\]
the limit $z\to1$ for the second term is 
\begin{align*}
\di\lim_{z\to 1}{
\k{-(z-1)\cot(\pi z)\ff{\phi'(z)}{\phi(z)^{2}}}
}
&=-
\di\lim_{w\to 0}
{w\cot(\pi w) \ff{\phi'(w-1)}{\phi(w-1)^{2}}
}
	\\&=-\ff{1}{\pi}\ff{m(2m-1)}{(2m)^{2}}
	=-\ff{2m-1}{4m\pi}.
\end{align*}

\end{proof}

\subsection{Proof of main theorem 2}
\renewcommand{\g}{\gamma}

Let $\b=\exp\k{\ff{\pii}{m}}$ as above.
\begin{thm}\label{th2}
\[
\z_{-1}(2m)=
\f{\,1\,}{2}+\ff{2m-1}{4m}-
\ff{\pi}{4m}
\sum_{1\le k\le 2m-1, k\ne m}
\b^{k}\cot(\pi \b^{k}).
\]
\end{thm}


\begin{proof}
Let $N, C(N), D(N)$ be as in the previous section. Again, the residue theorem claims that 
\[
\di \int_{C(N)}{g(z)}\,dz=
2\pi i \sum_{
\substack{\te{$a$:pole of $g$}\\
a\in D(N)}
 }\res(g, a).
\]
Taking the limit $N\to \mg$, the integral converges to 0 likewise.
It follows from Lemma \ref{l3} that 
\begin{align*}
	0&=2\pi i \sum_{a:\te{pole of $g$}}\res(g, a),
	\\0&=
	\res(g, 0)+\res(g, 1)+\res(g, -1)+
	\dsum_{n=2}^{\mug} \k{\res(g, n)+ \res(g, -n)}
	\\&+
\sum_{1\le k\le 2m-1, k\ne m}
\ff{\p}{2m}
\b^{k}\cot(\pi \b^{k})
	\\&=
	-\ff{1}{\p}-\ff{2m-1}{2m}\pi
	+\ff{2}{\pi}\z_{-1}(2m)+
\sum_{1\le k\le 2m-1, k\ne m}
\ff{\pi}{2m}
\b^{k}\cot(\pi \b^{k}).	
\end{align*}
Conclude that 
\[
\z_{-1}(2m)=
\f{\,1\,}{2}+\ff{2m-1}{4m}-
\ff{\pi}{4m}
\sum_{1\le k\le 2m-1, k\ne m}
\b^{k}\cot(\pi \b^{k}).
\]
\end{proof}

\subsection{Example}

\begin{ex}
For 
$s=4, m=2, \b=\exp\k{\tf{\pi i}{2}}=i$, 
an expression of real trigonometric function for $\z_{-1}(4)$ is 
\[
\z_{-1}(4)=
\ff{1}{2}+\ff{3}{8}-\ff{\pi}{8}
\k{i\cot(\pii)+i^{3}\cot(\pi i^{3})}
=
\ff{7}{8}-\ff{\pi}{4}\coth(\pi).
\]
Notice that this also shows that 
\[
\dsum_{k=1}^{\mug}\k{\z(4k)-1}
=\ff{7}{8}-\ff{\pi}{4}\coth(\pi)
\]
as in \cite[p.263]{bbc}. 
For $s=6, m=3, \b=\exp\k{\tf{\pi i}{3}}$, 
we see that 
\[
\z_{-1}(6)=\ff{1}{2}+\ff{5}{12}\]
\[-\ff{\pi}{12}
\k{
e^{\pi i/3}\cot(\pi e^{\pii/3})
+
e^{-\pi i/3}\cot(\pi e^{-\pii/3}
)
+
e^{2\pi i/3}\cot(\pi e^{2\pii/3})
+
e^{-2\pi i/3}\cot(\pi e^{-2\pii/3}
)}
\]
\[=
\ff{11}{12}-\ff{\p}{12}
\k{
\ff{\rt{3}\sinh \rt{3}\p}{\cosh\rt{3}\pi+1}
+
\ff{\rt{3}\sinh \rt{3}\p}{\cosh\rt{3}\pi+1}
}
\]
\[=
\ff{11}{12}-\ff{\p}{12}
\k{
\,2\rt{3}\,
\ff{2\sinh \tff{\rt{3}}{2}\p\cosh \tff{\rt{3}}{2}\p}
{2\cosh^{2} \tff{\rt{3}}{2}\p}
}
\]
\[
=\ff{11}{12}-\ff{\rt{3}}{6}\pi\tanh{\ff{\rt{3}}{2}\p}.
\]

\end{ex}

\paragraph*{\bf Acknowledgment:}
The first author thanks Satomi Abe, Shoko Asami and Michihito Tobe for sincerely supporting him. 
Also, the authors would like to thank Shigeru Iitaka. 
This research project arose from his online seminar in  2020-2022. 

\paragraph*{\bf Author Contributions:}
All authors contributed equally to the writing of this paper. All authors read and approved the final manuscript.
\paragraph*{\bf Conflicts of Interest:}``The authors declare no conflict of interest."

\end{document}